\documentclass{jlts}

\usepackage{amsmath,amssymb,amssymb,color,epsfig}
\usepackage[dvipsnames]{xcolor}

\title{SHORT NOTE ON AN OPEN PROBLEM}                                     
\author{Mohamed Bouali }                 
\lastname{Bouali}  



\address{Author name: Mohamed Bouali\\       \\        
Author address: Institut Pr{\'e}paratoire aux {\'e}tudes d'ingenieurs de Tunis.\\
  Facult\'e des Scinces de Tunis,  Campus Universitaire El-Manar, 2092 El Manar Tunis.\\\\            
Author email: bouali25@laposte.net                
}

\begin{document}


\maketitle
\begin{abstract}
In this work, we investigate a problem posed by Feng Qi and Bai-Ni Guo in their paper Complete monotonicities of functions involving the gamma and digamma functions.
\end{abstract}
\section{Introduction and statement of the main results}
A function $f$ is said to be {\it completely monotonic} on an interval $I$ if $f$ has derivatives
of all orders on $I$ which alternate successively in sign, that is
$$(-1)^nf^{(n)}(x)\geq 0,$$
for $x\in I$ and $n\geq 0$. If the inequality above is strict for all $x\in I$ and for all $n \geq 0$, then $f$
is said to be strictly completely monotonic. For a good monograph about this subject see for instance, \cite{sch} and \cite{wi}. In \cite{qi}, the author study a class of monotonic function involving the gamma and digamma functions. He shows  that the function $f_\alpha(x)=(x+1)^\alpha/\big((\Gamma(x+1))^{1/x}\big)$ is strictly completely monotonic on $(-1,+\infty)$ provided that $\alpha\leq 1/(1+\tau_0)<1$, where $\tau_0$ is the maximum of the function defined on $\N\times[0,+\infty)$ by
$$\tau(n,t)=\frac1n\Big(t-(t+n+1)\big(\frac1{1+t}\big)^{n+1}\Big).$$
 For $n=2,3$, the author compute numerically the maximum, his finds $\max\limits_{t>0}\tau(2,t)\simeq0.264076$, and
 $\max\limits_{t>0}\tau(3,t)\simeq 0.271807$. At my knowledge, no analytic method is given to compute $\tau_0$. 

 He posed the problem: find the maximum value
 $$\alpha=\max_{(n,t)\in\N\times(0,+\infty)}\tau(n,t).$$
 Using mathematica computations the author conjectured in \cite{qi} that $\alpha>0.298$.

  In this work, we answered the question, and we show that $\alpha$ is given by
  $$\alpha=\frac{\ell}{1+\ell+\ell^2},$$
  where $\ell$ is the unique solution of following equation
  $$e^{\frac1\ell}-\frac1{\ell^2}-\frac1\ell-1=0.$$
  Numerically $\ell\simeq 0.5576367386...$, and $\alpha=0.2984256075...$

  First of all, we show for every $n\geq1$, there is a unique sequence $t_n\in(0,n)$ which is increasing and such that
  $\max\limits_{t>0}\tau(n,t)=\tau(n,t_n):=\alpha_n$, and in second time, we prove that the sequence $\alpha_n$
converges towards growth to $\alpha$. Which allows us to deduce that $$\max_{(n,t)\in\N\times(0,+\infty)}\tau(n,t)=\alpha.$$
\section{Proof of main results}
Let's define on $[1,+\infty)\times[0,+\infty)$ the function
$$\tau(x,t)=\frac1x(t-(t+x+1)(\frac t{1+t})^{x+1}).$$
\begin{proposition}\label{0p}  For every $x\geq 1$, the function $t\mapsto\tau(x,t)$ attains its maximum in only one point $t(x)\in(0,x)$. The value of the maximum $\alpha(x):=\tau(x,t(x))$ is given by
$$\alpha(x)=\frac{(1+x) t(x)}{x^2+(1+t(x))^2+x (2+t(x))}$$
\end{proposition}

{\bf Proof.} 1)
Deriving the function $\tau(x,t)$ with respect to $t$, gives
$$\partial_t\tau(x,t)=\frac1x(1-(\frac t{1+t})^{x+1}-(t+x+1)(x+1)(\frac1t-\frac1{t+1})(\frac t{1+t})^{x+1}),$$
$$\partial_t\tau(x,t)=\frac1x(1-(\frac t{1+t})^{x+1}-(t+x+1)(x+1)\frac1{t(1+t)}(\frac t{1+t})^{x+1}).$$
The second derivative is given as follow
$$\partial_t^2\tau(x,t)=(1+x) t^{-1+x} (1+t)^{-3-x} (-1-x+t).$$
Hence, $\partial_t^2\tau(x,t)=0$ on $(0,+\infty)$ if and only if $t=x+1$. Moreover, $\partial_t\tau(x,0)=\frac1x$, and
$$\partial_t\tau(x,x)=\frac1x(1-3(\frac x{1+x})^{x+1}-\frac1{x}(\frac x{1+x})^{x+1}). $$
For $x\geq 1$, consider the function $h(x)=1-(3+\frac1x)(\frac x{x+1})^{x+1}$. We have by successive differentiation
$$h'(x)=-x^x (1 + x)^{-1 - x}h_1(x) ,$$
where $h_1(x)=3 + (1 + 3 x) \log(1 + \frac1x)$,
$$h_1'(x)=\frac{1}{x}+\frac{2}{1+x}+3 \log(x)-3 \log(1+x),$$
and $$h_1''(x)=\frac{  x-1}{x^2 (1 + x)^2}\geq 0.$$
Since, $\lim_{x\to+\infty}h_1'(x)=0,$ hence $h_1'(x)\leq0$. Moreover,
$\lim_{x\to+\infty}h_1(x)=0$, and hence $h_1(x)\geq 0$. One deduces that $h$ is non increasing on $[1,+\infty[$, and $h(1)=0$.
Thus, for every $x> 1$ $$\partial_t\tau(x,x)<0.$$
 Which implies that the derivative $t\mapsto\partial_t\tau(x,t)$ decreases strictly on $(0,x+1)$. By the fact that the second derivative acroses the $x$-axis only one time on $(0,+\infty)$, then there is a unique $t(x)\in(0,x+1)$ such that $$\partial_t\tau(x,t(x))=0,$$
 and $$t(x)\in(0,x),$$ in view of $\partial_t\tau(x,x)<0.$

Remark that in the interval $(0,x)$, $\partial_t^2(x,t)<0,$ hence the critical point $t(x)$ is a maximum for the function $t\mapsto\tau(x,t)$.

2) Recall that
$$\alpha(x)=\tau(x,t(x)).$$
 Furthermore, \begin{equation}\label{e}1-(\frac {t(x)}{1+t(x)})^{x+1}-(t(x)+x+1)(x+1)\frac1{t(x)(1+t(x))}(\frac {t(x)}{1+t(x)})^{x+1}=0.\end{equation}
and
\begin{equation}\label{eei}\alpha(x)=\frac1x(t(x)-(t(x)+x+1)(\frac {t(x)}{1+t(x)})^{x+1}).\end{equation}
 By substituting equation \eqref{e} in \eqref{eei}, we get
 \begin{equation}\label{ee}\alpha(x)=\frac{(1+x) t(x)}{x^2+(1+t(x))^2+x (2+t(x))}.\end{equation}
 \begin{lemma}\label{lill} The function $x\mapsto t(x)$ is well defined, 
  defines a $C^1$-diffeomorphism and satisfies
 $$\frac{(x+1)^2}{2x+3}\leq t(x)<x.$$
 \end{lemma}
 {\bf Proof.} 1) $x\mapsto t(x)$ defines a function by uniqueness proved in Proposition \ref{0p}.

 2) Let $a>0$, we saw that $\partial_t\tau(a,t(a))=0$, and $\partial^2_t\tau(a,t(a))<0$. Applying implicit theorem, then there is a neighborhood $V_a$ of $a$ and neighborhood $W_{t(a)}$ of $t(a)$ and $C^1$-diffeomorphism $\varphi:V_a\rightarrow W_{t(a)}$ such that for every $x\in V_a$, and $y\in W_{t(a)}$
 $$\partial_t\tau(x,y)=0\Leftrightarrow y=\varphi(x).$$
 By uniqueness $\varphi(x)=t(x)$.

 3) The right inequality $t(x)<x$ has been proved. To show the left one, it is enough to prove that
 $$\partial_t(x,\frac{(x+1)^2}{2x+3})\geq 0,$$
 in view of the decay of $t\mapsto\partial_t\tau(x,t)$ throughout the interval $(0,t(x))$, and the definition of $t(x)$, and uniqueness the result follows.

 Set $$\phi(x)=\partial_t(x,\frac{(x+1)^2}{2x+3}).$$
 By some algebra we get
 $$\phi(x)=\frac{(x+1)^{2x+2}}{x(x+2)^{2x+4}}\Big(\frac{(x+2)^{2x+4}}{(x+1)^{2x+2}}-7x^2-21x-16\Big).$$
 Let $$\Phi(x)=2\log(x+2)+2(x+1)\log(\frac{x+2}{x+1})-\log(7x^2+21x+16).$$
 Differentiate gives
 $$\Phi'(x)=2\log(\frac{x+2}{x+1})-\frac{21+14 x}{16+21 x+7 x^2},$$
 and $$\Phi'' (x)=\frac{-35 x^2 -105x-78}{(1+x) (2+x) (16+7 x (3+x))^2}<0.$$
 Since $\lim_{x\to+\infty}\Phi'(x)=0$, hence $\Phi'(x)>0$ for every $x\geq 1$. Thus,
 $$\Phi(x)\geq\Phi(1)>\Phi(0)=4\log2-\log16=0.$$
 Which gives, $$\phi(x)>0,$$
 and the result follows.
 \begin{proposition}\label{pip} The function $x\mapsto\alpha(x)$ is of class $C^1$ non decreasing on $[1,+\infty[$, with derivative
  $$\alpha'(x)=\partial_1\tau(x,t(x)),$$
  and satisfies for every $x\geq 1$ the inequality
  $$0\leq\alpha(x)<\frac{x}{3x+1}.$$
 \end{proposition}
 {\bf Proof.} 1) Recall that
$$\alpha(x)=\tau(x,t(x)).$$
$\alpha(x)$ is $C^1$ as composed of the $C^1$ functions $x\mapsto(x,t(x))$ and $(x,t)\mapsto f(x,t)$.

By differentiation we get
$$\alpha'(x)=\partial_1\tau(x,t(x))+t'(x)\partial_2\tau(x,t(x))=\partial_1\tau(x,t(x)),$$
where we used $\partial_2\tau(x,t(x))=0$. Deriving the expression of the function $\tau(x,t)$ with respect to $x$, gives

$$\partial_1\tau(x,t(x))=\frac{t(x) \left(-1-t(x)+\left(\frac{t(x)}{1+t(x)}\right)^x \big(1+t(x)+x (1+t(x)+x) \log(1+\frac1{t(x)})\big)\right)}{(1+t(x)) x^2}$$
$$=(1+t(x))\Big(-1+\frac{(1+t(x))(1+t(x)+x (1+t(x)+x) \log(1+\frac1{t(x)}))}{t(x)(1+t(x))+(t(x)+x+1)(x+1)}\Big),$$
where we used equation \eqref{e}.

For $0<u\leq x$, and $x\geq 1$, let's define 
 \begin{equation}\label{ii}k(u)=-1+\frac{(1+u)(1+u+x (1+u+x) \log(1+\frac1{u}))}{u(1+u)+(u+x+1)(x+1)}.\end{equation}
 First of all $$\alpha'(x)=\partial_1\tau(x,t(x))=(1+t(x))k(t(x)).$$
So, it is enough to show that $k(t(x))\geq 0$.

Differentiate yields
$$k'(u)=\frac{x \Big(u^2-(1+x)^3-u x (3+2 x)+u x (1+x) (2+2 u+x) \log\big(1+\frac{1}{u}\big)\Big)}{u \big(u^2+(1+x)^2+u (2+x)\big)^2}.$$
   Using the inequality $\log(1+1/u)\leq 1/u$, we get
   $$k'(u)\leq \frac{(1+u) (u-1-x) x}{u \left(u^2+(1+x)^2+u (2+x)\right)^2}<0.$$
   Thus $k(u)$ decreases, since $t(x)<x$, hence for every $x\geq 1$
   $$k(t(x))\geq k(x)=-1+\frac{(1+x)(1+x+x (1+x+x) \log(1+\frac1{x}))}{x(1+x)+(x+x+1)(x+1)}.$$
   Which gives

  \begin{equation}\label{ie}k(t(x))\geq -1+\frac{1+x+x(1+2x) \log(1+\frac1{x})}{3x+1}.\end{equation}
 Set $$\Theta(x)=1+x+x(1+2x) \log(1+\frac1{x})-3x-1,$$
 Differentiate, by straightforward computation it yields
 $$\Theta'(x)=-4+\frac{1}{1+x}+(1+4 x) \log\left(1+\frac{1}{x}\right),$$
 $$\Theta''(x)=\frac{-1-2 x (3+2 x)}{x (1+x)^2}+4 \log\left(1+\frac{1}{x}\right),$$
 and
 $$\Theta'''(x)=\frac{1-x}{x^2 (1+x)^3}.$$
 For $x\geq 1$, $\Theta'''(x)<0$, and $\Theta''(x)\geq \lim_{x\to+\infty}\Theta''(x)=0$. Hence
 $\Theta'(x)\leq\lim_{x\to+\infty}\Theta'(x)=0$. Then $\Theta$ decreases, moreover, $\Theta(x)=1/(6x)+o(1/x)$. then
 $\Theta(x)\geq 0$.
  By equation \eqref{ie}, one deduces the positivity of $k(t(x))$, namely for every $x\geq 1$, $k(t(x))\geq 0$. and,
   $$\alpha'(x)=\partial_1\tau(x,t(x))\geq 0.$$
2) We saw by Proposition \ref{0p}, $\alpha(x)=\psi(t(x)),$
where $$\psi(u)=\frac{(x+1)u}{x^2+(u+1)^2+x(u+2)}.$$
Deriving with respect to $u$, it yields
 $$\psi'(u)=\frac{(1+n)^2-u^2}{\big(n^2+(1+u)^2+n (2+u)\big)^2}\geq 0,$$
 for every $u\geq x+1$.

 Since, by Lemma \ref{lill}, $t(x)< x$, and $\psi(x)=x/(3x+1)$. One deduces that, for every $x\geq 1$
 $$0\leq\alpha(x)\leq\frac{x}{3x+1}.$$

 \begin{proposition}
 \begin{description}
 \item[(i)] The sequence $t_n$ increases, and the sequence $\frac{t_n}n$ converges to $\ell$, the unique solution of the equation
$$e^{\frac1\ell}-\frac1{\ell^2}-\frac1\ell-1=0.$$
Numerical computation gives $\ell\simeq0.5577$.

 \item[(ii)] The sequence $\alpha_n$ is bounded. Moreover, $\alpha_n$ converges to $\alpha$, where
 $$\alpha=\frac{\ell}{1+\ell+\ell^2}.$$
 Numerically $\alpha\simeq 0.298438.$
 \item[(iii)] $$\max_{(n,t)\in\N\times(0,+\infty)}\tau(n,t)=\alpha.$$
 \end{description}
 \end{proposition}
 {\bf Proof.} i)a) Let $t_{n+1}$ denotes the unique zero of the function $t\mapsto\partial_t\tau(n+1,t)$.
 Straightforward computation gives $$\partial_t\tau(n,t_{n+1})=\frac1n\Big(1-\big(1+\frac{(t_{n+1}+n+1)(n+1)}{{t_{n+1}}(1+t_{n+1})}\big)(\frac {t_{n+1}}{1+t_{n+1}})^{n+1}\Big).$$
 Furthermore, by equation $\partial_t\tau(n+1,t_{n+1})=0,$ we get
 $$(\frac {t_{n+1}}{1+t_{n+1}})^{n+2}=\frac{t_{n+1}(1+t_{n+1})}{t_{n+1}(1+t_{n+1})+(t_{n+1}+n+2)(n+2)},$$
 Which implies that
 $$\partial_t\tau(n,t_{n+1})=
 \frac1n(1-\frac{t_{n+1}(t_{n+1}+1)+(t_{n+1}+n+1)(n+1)}{t_{n+1}(t_{n+1}+1)+(t_{n+1}+n+2)(n+2)}\frac {t_{n+1}+1}{t_{n+1}}).$$
 In other words, $$\partial_t\tau(n,t_{n+1})=-\frac{(1+n) (1+n-t_{n+1})}{nt_{n+1} \left((2+n)^2+(3+n) t_{n+1}+t^2_{n+1}\right)}< 0,$$
 in view of $t_{n+1}\leq n+1$.

 Using the fact that $t\mapsto\partial_t\tau(n,t)$ decreases on $(0,n+1)$, $t_n,t_{n+1}\in(0,n+1)$ and  $\partial_t\tau(n,t_n)=0$, we get
 $t_n< t_{n+1}$

 b)
   First of all, by Lemma \ref{lill}, $(n+1)^2/n(2n+3)\leq t_n/n\leq1$ and the left hand side is bounded, then $t_n/n$ is bounded too. Let $t_{n_k}/n_k$ be a some convergent subsequence, and $\ell=\lim_{k\to\infty}t_{n_k}/n_k\in[1/2;1]$.

Using equation \eqref{e} with $x=n_k$, one deduces that, as $k\to+\infty$
$$\frac{e^{-1/\ell} (-1-\ell+(-1+e^{\frac{1}{\ell}}) \ell^2)}{\ell^2}=0,$$
or if we set $a=\frac1\ell$, with $a\geq 1$,
\begin{equation}\label{j}e^{a}=a^2+a+1.\end{equation}

Let $\eta(a)=e^a-a^2-a-1$, then $\eta'(a)=e^a-2a-1$, and $\eta''(a)=e^a-2>0$ for $a\geq1$. Hence, $\eta'$ increases on $(1,+\infty)$. Since, $\eta'(1)=e-3<0$, and $\eta'(2)=e^2-5>0$. Hence there is a unique $a_0\in]1,2[$, $\eta'(a_0)=0$, $a_0\sim1.26$.
Remark that $\eta(1)=e-3<0$, and $\eta(a_0)=2a_0+1-a_0^2-a_0-1=a_0(1-a_0)<0$. Moreover on $(1,a_0]$ the function $\eta$ decreases and is strictly negative, and increases on $(a_0,+\infty)$ with $\lim_{a\to+\infty)}\eta(a)=+\infty$. Thus equation \eqref{j} admit a unique solution $x_0$ in $(a_0,+\infty)$. Which implies that $\ell$ is the unique limit of a subsequence and the sequence $t_n/n$ converges to $\ell$. Numerically $x_0\simeq1.793$, and $\ell=1/x_0\simeq 0.5577,$



b) Form Proposition \ref{0p}, with $x=n$, $\alpha_n=\alpha(n)$, and $t_n=t(n)$, $$\alpha_n=\frac{(1+n) t_n}{n^2+(1+t_n)^2+n (2+t_n)},$$

As $n$ goes to $+\infty$, and the fact that $t_n/n\to\ell$, one gets $$\alpha:=\lim_{n\to\infty}\alpha_n=\frac{\ell}{1+\ell+\ell^2}.$$
 Numerically $\alpha\simeq0.298438.$

iii) We saw by Proposition \ref{pip} that the sequence $\alpha_n$ increases and converges to $\alpha=\sup_n\alpha_n$. First of all, for every $n\geq 1$,
$$\tau(n,t)\leq\max_{t>0}\tau(n,t)=\alpha_n\leq\alpha.$$
Hence, $\max_{(n, t)\in\N\times(0,+\infty)}\tau(n,t)$ is well defined and
 $$\max_{(n, t)\in\N\times(0,+\infty)}\tau(n,t)\leq\alpha. $$
Moreover,
$$\alpha_n=\tau(n,t_n)\leq\max_{(n, t)\in\N\times(0,+\infty)}\tau(n,t).$$
One deduces that
$$\max_{(n, t)\in\N\times(0,+\infty)}\tau(n,t)=\sup_{n\geq 1}\alpha_n=\alpha.$$
\begin{figure}[h]
\centering\scalebox{0.6}{\includegraphics[width=15cm, height=10cm]{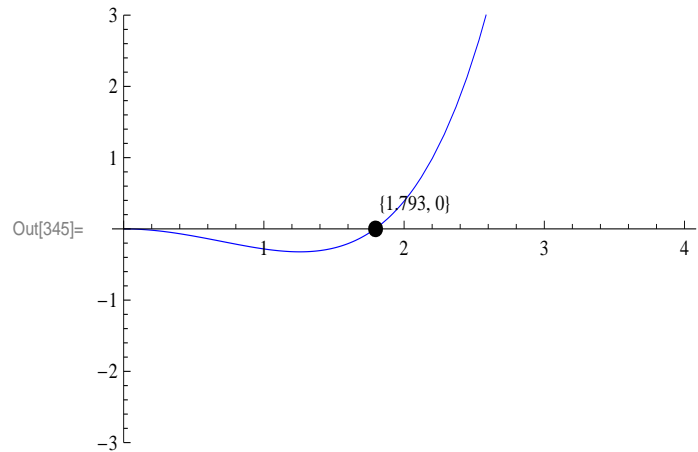}}

\caption{\bf Plot of the function $f(a)=e^a-a^2-a-1$, and $x_0=1.793$.}
\end{figure}

\begin{figure}[h]
\centering\scalebox{0.6}{\includegraphics[width=15cm, height=10cm]{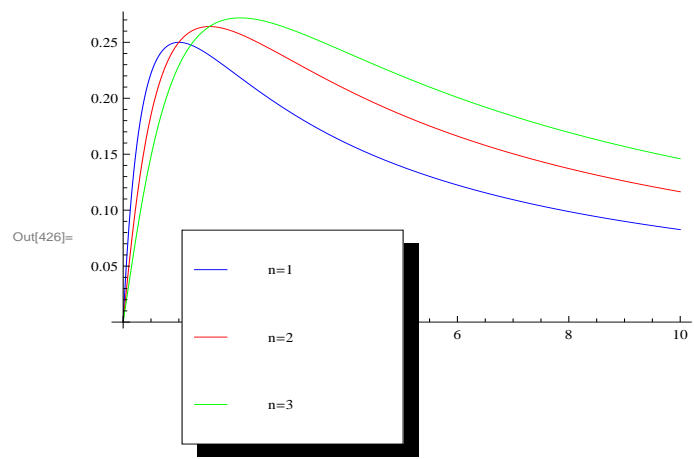}}

\caption{\bf Plot of the function $\tau(n,t)$, $n=1,2,3$}
\end{figure}

\end{document}